\documentclass[12pt]{amsart}
\usepackage{amsmath,amssymb,amsthm}
\title[a. c. subspace of a nonselfadjoint
operator]{Notions of absolutely continuous subspace for
nonselfadjoint operators}
\author{Roman Romanov}
\email{morovom@gmail.com}
\address{Laboratory of Quantum Networks,
        Institute for Physics,
        Saint Peters\-burg State University,
        198504 St. Petersburg
        RUSSIA.}

\begin{document}

\begin{abstract}
We give an example of an operator with different weak and strong
absolutely continuous subspaces, and a counterexample to the
duality problem for the spectral components. Both examples are
optimal in the scale of compact operators.

\end{abstract}

\maketitle

\newtheorem{theorem}{Theorem}
\newtheorem{proposition}[theorem]{Proposition}
\newtheorem{corollary}[theorem]{Corollary}
\newtheorem{lemma}[theorem]{Lemma}
\def\thelemma{\unskip}

\newtheorem{remark}{Remark}
\def\theremark{\unskip}
\newtheorem{definition}{Definition}

\newcommand{\C}{\mathbb{C}}
\newcommand{\D}{\mathbb{D}}
\newcommand{\R}{\mathbb{R}}
\newcommand{\N}{\mathbb{N}}
\newcommand{\Z}{\mathbb{Z}}

\def\be{\begin{equation}}
\def\ee{\end{equation}}
\def\bequnan{\begin{eqnarray*}}
\def\eequnan{\end{eqnarray*}}
\def\la{\label}
\def\von{\varepsilon}
\def\hp{\widehat{p}}
\def\vsup{\operatorname{ess\ sup}}
\def\vinf{\operatorname{ess\, inf}}
\def\dist{\operatorname{dist}}
\newcommand{\bS}{{\bf S}}
\def\cal{\mathcal}
\def\cD{{\cal D}}
\def\cP{{\cal P}}
\def\cR{{\cal R}}
\def\cN{{\cal N}}
\def\cM{{\cal M}}
\def\cH{{\cal H}}
\def\cL{{\cal L}}
\def\cK{{\cal K}}
\def\cA{{\cal A}}
\def\tg{\tilde{g}}
\def\tu{\tilde{u}}
\def\tv{\tilde{v}}
\def\wt{\widetilde}
\def\tL{\widetilde{L}}
\def\tV{{\widetilde{V}}}
\def\tm{\widetilde{m}}
\def\wh{\widehat}
\def\Ran{\operatorname{Ran}}
\def\slim{\operatorname{s-lim}}
\def\({\left(}
\def\){\right)}
\def\H{{\bf H}}
\def\llangle{\left\langle}
\def\rrangle{\right\rangle}
\def\len{\left\|}
\def\rin{\right\|}
\def\im{\mathrm {Im}\,}
\def\re{\mathrm {Re}\,}

{\bf Introduction}. The known definitions of the absolutely
continuous (a. c.) subspace for nonselfadjoint operators fall into
two groups: the weak ones and the strong one. The strong
definition was first suggested by L. Sakhnovich \cite{Sachn} in
the case of dissipative operators.

\begin{definition}
Let $ L $ be a completely non-selfadjoint dissipative operator in
a Hilbert space $ H $ with a bounded imaginary part $ V $. Any of
the following coinciding subspaces is called the strong a. c.
subspace $ H_{ac} $ of $ L $,
\begin{equation*}
\begin{array}{cl}
 (i) & \text{The invariant subspace of $ L $ corresponding to
the canonical } \\ & \text{factorization of its characteristic
function;} \\ (ii) & \text{The minimal subspace containing all the
invariant subspaces $ X $} \\ & \text{of $ L $ such that $ \left.
L \right|_X = W A W^{ -1 } $ for an a. c. selfadjoint operator $ A
$} \\ & \text{and a bounded and boundedly invertible operator $ W
\colon X \to X $;} \\(iii) & \mathop{Clos } \left\{ u \in H \colon
\; \left. V^{ 1/2 } \( L - z \)^{ -1 } u \right|_{ \C_+} \in
\H^2_+ \right\}.
\end{array}
\end{equation*}
\end{definition}

The notation we use is given in the end of Introduction. The
corres\-ponden\-ce meant in definition (i) is the one between
invariant subspaces of an operator and regular factorizations of
its characteris\-tic function in the framework of the
Sz\"{o}kefalvi - Nagy - Foia\c{s} functional model. The
equivalence of (ii) and (iii) is a corollary of this
corres\-ponden\-ce. The equivalent definition (iii) in model-free
terms was found in \cite{N}. In this latter form, the definition
was generalized to non-dissipative perturbations of selfadjoint
operators \cite{N} and non-contractive perturba\-tions of unitary
ones \cite{MV}.

The weak definitions of the a. c. subspace are obtained if we try
to generalize directly the "selfadjoint" definition using the
property expressed by the Riesz brothers theorem as a substitute
for the absolute continuity of (non-existent, in general) spectral
measure. This leads to spaces defined by the requirement that the
matrix element of the resolvent be of the Hardy class $ H^p $.
They were first introduced and studied in \cite{Tikhon,Tikhang}
and are called weak a. c. subspaces. A weak a. c. subspace
contains the strong one because the restriction of the operator to
the a. c. subspace is quasi-similar to an a. c. selfadjoint
operator. A natural question is whether these subspaces coincide.
So far, it was answered in affirmative in two situations ($ p =
2)$:

(\textit{i}) When $ L $ is dissipative  \cite{equival};

(\textit{ii}) When the characteristic function of $ L $ has weak
bounda\-ry values a. e. on the real axis \cite{Ryzov}. This holds
true, for instance, for trace class perturbations of a selfadjoint
operator, and, more generally, if the function admits scalar
multiple.

The first main result of the present paper given by theorems
\ref{example}  and \ref{exampleline} is ---
\begin{itemize}
\item An example of an operator with different weak and strong
a. c. subspaces.
\end{itemize}

The example in theorem \ref{exampleline} is a (non-dissipative)
perturbation of a self\-ad\-joint operator. Theorem \ref{example}
provides an analogous result for perturbations of unitary
operators. The operators we construct are in fact similar to the
standard bilateral shift \cite{Sh} in the unit circle case and to
the generator of bilateral shift on $ \mathbb{R} $ in the real
line case.

Whichever definition of the a. c. subspace is used, it is natural
to ask whether the orthogonal complement of it coincides with the
singular subspace of the adjoint operator. The latter is defined
to be the closure of the set of vectors such that the matrix
element of the resolvent on such a vector has zero jump a. e.
while crossing the essential spectrum (real line and unit circle
in the case of perturbations of selfadjoint and unitary operators,
respectively). This question is sometimes referred to as the
duality problem for spectral components and is known to have
affirmative answer in the case of trace class perturbations
\cite{MV,Tikhon}. Our second main result is theorem
\ref{dualexample} ---

\begin{itemize}
\item An example of an operator with trivial singular subspace of the
adjoint operator and nontrivial orthogonal complement of the weak
a. c. subspace.
\end{itemize}

The example can be described as a non-con\-tracti\-ve conjunction
of two bilateral shifts.

The examples in theorems \ref{example} and \ref{dualexample} are
optimal in the sense that they are additive perturbations of
unitary operators by an operator arbitrarily close to the trace
class in the sense that its $s$-numbers can be chosen to be
estimated above by an arbitrary given monotone non-summable
sequence.

Theorems \ref{example}  and \ref{exampleline}  are proved in \S 1,
theorem \ref{dualexample} --- in \S 2. In \S 3 we analyze the weak
definition of the a. c. subspace for $ p \ne 2 $, and show that in
the dissipative case it gives the same subspace as for $ p = 2 $.

NOTATION:

$ \mathbb{C}_\pm = \{ z : \pm \im z > 0 \} $, $ \mathbb{T} = \{ z
\colon | z | = 1 \} $; $ \mathbb{D} = \{ z \colon |z| < 1 \} $.

$ H^p_\pm $, $ 0 < p \le 2 $, - Hardy classes of analytic
functions in $ \C_\pm $, respectively.

$ H^2 $ - the Hardy space for the unit disk.

For functions with values in a Hilbert space $ H $,

$ \H^2_\pm $ - the Hardy space of $ H $-valued functions in $
\C_\pm $, respectively. The norm of an $ f \in \H^2_\pm $ is given
by $\sup_{ \von > 0 } \int_\R \len f ( k \pm i\von ) \rin_H^2 dk <
\infty $;

$ \H^2 $ - the Hardy space of $H$-valued functions in the unit
disk.

$ \{ e_n \} $ - the standard basis in $ l^2 ( \mathbb{Z} ) $.

For an operator $ L $ in a Hilbert space,

$ \cD ( L ) $ - the domain of $ L $;

$ f_{ u , v } ( z ) = \llangle \( L - z \)^{ -1 } u , v \rrangle
$.

$ \mathbf{S}^p $, $ p \ge 1 $, - the Shatten - von Neumann classes
of compact operators with summable $ p $-th power of their
singular numbers.

The subscripts $ _\pm $ with functions in the complex plane stand
for their respective restrictions to $ \mathbb{C}_\pm $.

Various subspaces corresponding to abstract operators are defined
in the paper. We will often suppress the explicit indication of
the operator in the notation for the subspaces when it is clear
which operator it refers to.

\bigskip

\S 1. Let $ L $ be a closed operator in a Hilbert space $ H $ such
that $ \sigma ( L ) \cap \C_\pm $ are discrete sets. For any $ p
\le 2 $ one can define the following subspaces in $ H $,
\begin{eqnarray} \label{defpla} & H_{ac}^{w, p } ( L ) \stackrel{ \mbox\small{def}}{ =
} \mathop{Clos} \wt{H_{ac}^{w,p}} ( L ) , & \\ \nonumber &
\wt{H_{ac}^{w,p}} ( L ) \stackrel{ \mbox\small{def}}{ = } \left\{
\begin{array}{cl} u \in H \colon & \( L - z \)^{ -1 } u \mbox{ is analytic in }
\C \setminus \R ,
\\ & \llangle \( L - z \)^{ -1 } u , v \rrangle_
\pm \in H^p_\pm \mbox{ for all } v \in H . \\
\end{array} \right\}.
\end{eqnarray}

In the case $ p = 2 $ this subspace is the \textit{weak a. c.
subspace} of the operator $ L $. Elements of the linear set $
\wt{H_{ac}^w} ( L ) $ are called {\it weak smooth vectors}. If $ L
$ is self-adjoint, then for all $ p, 1 < p \le 2 $, the subspaces
$ H_{ac}^{w, p } ( L ) $ coincide with the a. c. subspace of the
operator $ L $ defined in the standard way. We include a proof of
this folklore-type assertion in the Appendix. We shall omit the
index $ p $ in our notation in the case $ p = 2$ and write $
H_{ac}^w $ for $ H_{ac}^{w, 2 } $. It follows from the uniform
boundedness principle that for $ u \in \wt{H_{ac}^{w, p}} $ the $
H^p $-norms of the functions $ f_{ u , v } ( z ) $ are bounded
above when $ v $ ranges over the unit ball in $ H $.

For clarity, we restrict our consideration to the situation of the
perturbation theory. From now on, it is assumed that

\medskip

{\parindent=0cm {\bf (A)}} {\it $ L $ is a completely
nonself-adjoint operator of the form $ L = A + i V $, $ A = A^* $,
$ V = V^* $, $ \cD ( L ) \colon = \cD ( A ) \subset \cD ( V ) $,
and $ V $ is $ A $-bounded with a relative bound less than $ 1 $,
that is, $ \len V u \rin^2 \le a \len A u \rin^2 + b \len u\rin^2
$, $ a< 1 $, for all $ u \in \cD ( A ) $.}

\smallskip

This assumption implies, in particular, that \be\la{inft} i\tau \(
L + i \tau \)^{ -1 } \stackrel{s} \longrightarrow I, \; \tau \to
\pm \infty , \ee and that the operator $ \left| V \right|^{ 1/2 }
\( L - z \)^{ -1 } $ is defined as a bounded operator for all $ z
\in \rho ( L ) $.

\begin{definition}[\protect\cite{N}] \la{strongns}
The subspace \bequnan & H_{ac} ( L ) \stackrel{ \mbox\small{def}}{
= } \mbox{Clos } \wt{H_{ac}} ( L ) , & \\ & \wt{H_{ac}} ( L )
\stackrel{ \mbox\small{def}}{ = } \left\{
\begin{array}{cl} u \in H \colon & \( L - z \)^{ -1 } u \mbox{ is analytic in }
\C \setminus \R ,  \; \;
\\ & \(  \left| V \right|^{ 1/2 } \( L - z \)^{ -1 } u \)_\pm \in \H^2_\pm \\
\end{array} \right\}.
\eequnan is called the strong a. c. subspace of the operator $ L$.
Elements of the set $ \wt{H_{ac}} ( L ) $ are called strong smooth
vectors.
\end{definition}

Notice that there exists a natural generalization of this
definition applicable to operators which do not satisfy the
assumption (A) \cite{Ryzov}.

The main property of the strong smooth vectors is expressed by the
fol\-low\-ing

\medskip

{\parindent=0cm {\bf Proposition \cite[Theorem 4]{N}.}} {\it There
exists a Hilbert space $ \cN $, an a. c. self-adjoint operator $
A_0 $ in $ \cN $, and a bounded operator $ P \colon \cN \to H $
such that $ P \cN = \wt{H_{ac}} ( L ) $ and the equality $$ \( L -
z \)^{ -1 } P g = P \( A_0 - z \)^{ -1 } g. $$ holds for all $ g
\in \cN $ and $ z \notin \R $, $ z \in \rho ( L ) $.}

\begin{corollary} \la{inter}
$ H_{ac}^w ( L ) \supset H_{ac} ( L ) $.
\end{corollary}

A similar theory is available for perturbations of unitary
operators \cite{MV}. Let $ T $ be a bounded operator such that $
\sigma ( T ) $ has no accumulation points off $ \mathbb{T} $.

\begin{definition} \la{dewcon} The weak a. c.
subspace of the operator $ T $ is the set
\begin{eqnarray*} & H_{ac}^w ( T ) \stackrel{ \mbox\small{def}}{ =
} \mathop{Clos} \wt{H_{ac}^w} ( T ) , & \\ &  \wt{H_{ac}^w} ( T )
\stackrel{ \mbox\small{def}}{ = } \wt{H_+^w} ( T ) \cap \wt{H_-^w}
( T ),
\\ & \wt{H_+^w} ( T ) \stackrel{ \mbox\small{def}}{ = }
\left\{
\begin{array}{cl} u \in H \colon & \( T - z \)^{ -1 } u \mbox{ is analytic in }
\D , 
\\ & \left. \llangle \( T - z \)^{ -1 } u , v \rrangle \right|_{
\mathbb{D} } \in H^2 \mbox{ for all } v \in H \end{array}
\right\},
\\ & \wt{H_-^w} ( T ) \stackrel{ \mbox\small{def}}{ =
}\left\{\begin{array}{cl} u \in H \colon & \( T - z \)^{ -1 } u
\mbox{ is analytic in }
\C \setminus \overline{\D} , 
\\ &
\left. \llangle \( I - z T \)^{ -1 } u , v \rrangle \right|_{
\mathbb{D} } \in H^2 \mbox{ for all } v \in H
\end{array} \right\}.
\end{eqnarray*}

Let $ D_T \stackrel{ \mbox\small{def}}{ = } \left| I - T^* T
\right|^{ 1/2 } $. The subspace \bequnan & H_{ac}( T ) \stackrel{
\mbox\small{def}}{ = } \mbox{Clos } \wt{H_{ac}} ( T ) , & \\ &
\wt{H_{ac}} ( T ) \stackrel{ \mbox\small{def}}{ = } \left\{
\begin{array}{ccc}  u \in H: & (i) & \( T - z \)^{ -1 } u \mbox{ is analytic in }
\C \setminus \mathbb{T} , 
\\ & (ii) & \left. D_T \( T - z \)^{ -1 } u  \right|_{
\mathbb{D} } \in \H^2 , \hfill \\ & (iii) & \left. D_T \( I - z T
\)^{ -1 } u \right|_{ \mathbb{D} } \in \H^2 \hfill
\end{array} \right\}  \eequnan is called the (strong) a. c. subspace of the operator
$ T $. Elements of linear sets $ \wt{H_{ac}^w} ( T ) $ and $
\wt{H_{ac}} ( T ) $ are called weak and strong smooth vectors,
respectively.
\end{definition}

An analog of the proposition above holds for the strong smooth
vectors of $ T $ \cite{MV}.

\begin{corollary} \la{interc} $
H_{ac}^w ( T ) \supset H_{ac} ( T ) $.
\end{corollary}

We now proceed to construct the required examples, first for
perturbations of unitary operators. Given a selfadjoint operator $
D$, define $ \lambda_j ( D ) $ to be the eigenvalues of $ D $
enumerated in the modulus decreasing order. Let $ \{ \pi_n \} $, $
\pi_n > 0 $, be a monotone decreasing sequence.

\begin{theorem} \la{example}
There exists a bounded operator $ T $ obeying the following
con\-di\-tions,

\begin{itemize}
\item[\textit{i})]  $ H_{ac}^w ( T ) = H $,
\item[\textit{ii})] $ H_{ac} ( T ) = \{ 0 \} $,
\item[\textit{iii})] $ I - T^* T \in \mathbf{S}^p $ for all $ p > 1 $,
\item[\textit{iv})] $T$ is similar to a unitary operator.
\end{itemize}

Moreover, for any sequence $ \{ \pi_n \} \notin l^1 $ there exists
an operator $ T $ satisfying the conditions above with
\textit{iii}) replaced by
\begin{itemize}
\item[\textit{iii'})] $ \left| \lambda_n ( I - T^* T ) \right|
\le \pi_n $.
\end{itemize}
\end{theorem}

This theorem is optimal in the sense that the subspaces $ H_{ac}^w
$ and $ H_{ac} $ are known \cite[Proposition 4.10 and Theorem
C]{Tikhang} to coincide\footnote{The result in \cite{Tikhang} is
more general pertaining to operators with spectrum on a curve. In
the situation under consideration it is given in the unpublished
thesis \cite{Tikhon}.} if $ I - T^* T \in \mathbf{S}^1 $, provided
that $ \mathbb{D} \not\subset\sigma_{ess} ( T ) $. In the
terminology of \cite{GK}, the theorem says that no condition of
the form $ I - T^* T \in \mathbf{S}_\pi $ where $ \mathbf{S}_\pi $
is a symmetrically-normed ideal of compact operators containing $
\mathbf{S}^1 $ properly, guarantees the coincidence of $ H_{ac}^w
$ and $ H_{ac } $.

Notice that the similarity of operators respects weak a. c.
subspace, therefore the conditions \textit{i}) and \textit{iv})
combined are equivalent to saying that $ T $ is similar to an a.
c. unitary operator.

\begin{proof}
Let $ H = \ell^2 (\mathbb{Z}) $. We shall construct a sequence $
\left \{\rho_n \right\}_{ n= -\infty}^{ + \infty } $ of positive
numbers such that the weighted bilateral shift operator $ T $
defined by \[ T e_j = \rho_{ j-1 } e_{ j-1 } , j \in \mathbb{Z} ,
\]
has the required properties. Assume that $$  \sum_j \left| \rho_j
- 1 \right|^p < \infty \eqno(*) $$ for any $ p > 1 $. We are going
to need the explicit formula for the resolvent of $ T $, 
\[ \( \( T - \lambda \)^{ -1 } f\)_m = \left\{ \begin{array}{cc}
{\displaystyle \sum_{ k < m } f_k \frac{\lambda^{ m - k - 1 }}{
\prod_k^{ m-1 } \rho_j } \, ,} &  | \lambda | < 1 , \\
\noalign{\vskip5pt} {\displaystyle - \sum_{ k \ge m } f_k \frac{
\prod_m^{ k-1 } \rho_j }{\lambda^{ k - m + 1 }} \, ,} & | \lambda
| > 1 .\end{array} \right.
\] 

Proceeding, let us check that the following implication holds $$ \left. D_T \( T - z \)^{ -1 } u
\right|_{ \mathbb{D} } \in \H^2 \Longrightarrow \sum_{ n
> 0 } \frac{ \left| 1 - \rho_n^2 \right| } {\prod_0^{ n-1 } \rho_j^2 } <
\infty , \eqno(**) $$ provided that $ u \ne 0 $. This is done by a direct computation. In
the situation under considera\-tion
\[ D_T = \hbox{diag}\{ \left| 1 - \rho_n^2 \right|^{ 1/2 } \} . \]
We have: \bequnan \int_{ -\pi }^\pi \len D_T \( T - z \)^{ -1 } f
\rin^2 d\theta = \sum_n \left| 1 - \rho_n^2 \right| \int_{ -\pi
}^\pi \left| \sum_{ k < n } f_k \frac{ z^{ n - k - 1 }}{ \prod_k^{
n-1 } \rho_j } \right|^2 d\theta \\ = \sum_n \left| 1 - \rho_n^2
\right| \sum_{ k < n } r^{ 2(n - k - 1) } \frac{\left| f_k
\right|^2}{ \prod_k^{ n-1 } \rho_j^2} . \eequnan Thus, the
function $ D_T \( T - z \)^{ -1 }  \! f $ is in $ \H^2 $ if, and
only if, the quantity
\[ \sum_n \left| 1 - \rho_n^2 \right| \sum_{ k < n } \frac{\left|
f_k \right|^2}{ \prod_k^{ n-1 } \rho_j^2 } = \sum_k  \( \sum_{ n >
k } \frac{\left| 1 - \rho_n^2 \right|}{ \prod_k^{ n-1 } \rho_j^2 }
\) \left| f_k \right|^2 \] is finite. This means that the sum in
parentheses in the right hand side must be finite for some $ k $.
Since this sum, obviously, converges or diverges for all $ k $
simultaneously, the implication (**) is established.

The existence of an operator $ T $ enjoying the properties
\textit{i}) --- \textit{iv}) will be proved if we construct an
example of the sequence $ \{ \rho_j \} $ such that $ T $ is
similar to an a. c. unitary operator, the sum in (**) diverges,
and condition (*) is satisfied. Let $ a_j = 1 - 1/j $ and let
\bequnan& \rho_j = 1 , j \le 1, &
\\ & \rho_{2j} = a_j , j \ge 1, &  \\ & \rho_{2j+1} = a_j^{ -1 } , j \ge 1
. & \eequnan With this choice, (*) and the divergence of the sum
in (**) are straight\-for\-ward. Then, define \bequnan&
 w_{2j+1} = a_j^{ -1 } , j \ge 1 ; & \\
& w_j = 1  \; \, \text{ otherwise}. & \eequnan The diagonal
operator $ W = \mbox{diag} \{ w_j \} $, defined by the sequence $
\{ w_j \} $ in $ H $, is obviously bounded, boundedly invertible,
and it is easy to check that $ W^{ -1 } T W $ is the unitary
operator of (non-weighted) shift in $ H $.

To verify the second assertion of the theorem one can assume
without loss of generality that $ \pi_{2j + 1 } = \pi_{2j} $. It
is then enough to take  $ \{ a_j \} $ to be any sequence of
positive numbers such that $ a_j \to 1 $, $ | 1 - a_j | \le \pi_j
/2 $, but $ \sum | 1 - a_j | = \infty $, in the construction
above.
\end{proof}

\begin{remark} Theorem \ref{example} shows that the linear
resolvent growth condition \be\la{lrg} \sup_{ z \notin \mathbb{T}}
\( \big| 1 - |z| \big| \len \( T - z \)^{ -1 } \rin \) < \infty ,
\ee and even the combination of it and the property \textit{iii}),
do not imply the coinci\-dence of $ H_{ac} $ and $ H_{ac}^w $.
Also, it shows that in general simi\-lari\-ty of operators does
not respect the strong a. c. subspace.
\end{remark}

We now turn to perturbations of selfadjoint operators. In
principle, the corresponding example can be obtained by the Caley
transform. We prefer to give a direct construction based on the
same idea, because the Caley transform of the operator in theorem
\ref{example} does not belong to the class of operators for which
we have defined the strong a. c. subspace. This is not a real
hitch, however, for the Caley transform will indeed be the example
required if we use the general definition of $ H_{ac} ( L ) $ from
\cite{Ryzov} mentioned above.

Let $ H = L^2 ( \mathbb{R} ) $, and $ q ( x ) $ be a bounded
function on $\mathbb{R} $ satisfying\footnote{$ q ( x ) = \sin x /
x $ is the simplest example.}
\begin{eqnarray*} & \sum_n \( \int_n^{ n + 1 } | q |^2  \)^{ p / 2
} < \infty \mbox{ for all } p
> 1 , & \\ & q \notin L^1, & \\ & q \mbox{ is
conditionally integrable} . \end{eqnarray*}

Let $ L $ be the operator in $ H $ defined by the differential
expression \[  L =  i \frac d{dx}  + i q ( x ) \] on its natural
domain. Notice that the operator $ L $ is similar to the operator
$ A =  i \frac d{dx} $: \[ L = W A W^{ -1 } , \] where $ W $ is
the operator of multiplication by the function $ \exp \( - \int_{
-\infty }^x q \) $.

\begin{theorem} \la{exampleline}
The operator $ L $ obeys the following conditions,

\begin{itemize}
\item[\textit{i})]  $ H_{ac}^w ( L ) = H $,
\item[\textit{ii})] $ H_{ac} ( L ) = \{ 0 \} $,
\item[\textit{iii})] $ \( L - z \)^{ -1 } - \( A - z \)^{ -1 }
\in \mathbf{S}^p $ for all $ p > 1 $, $ \im z \ne 0 $,
\item[\textit{iv})] $L$ is similar to a selfadjoint operator.
\end{itemize}
\end{theorem}

\begin{proof} Since $ A $ is an a. c. selfadjoint operator,
\textit{i}) and  \textit{iv}) are immediate. Let $ V = \im L $.
Suppose that $ u $ is a strong smooth vector of $ L $. Since $ W $
commutes with the multiplication by a function, this is equivalent
to saying that the restrictions of the function \[  \varphi (z) =
\left| V \right|^{ 1/2 } \( A - z \)^{ -1 } g, \; g = W u , \]
belong to $ \H^2_\pm $ in the respective halfplanes. In turn, the
latter is equivalent to the condition \[ \int_{ \mathbb{R} } \len
\left| V \right|^{ 1/2 } e^{ itA} g \rin^2 dt < \infty \] by the
Parseval equality for the vector Fourier transform. We have
\[ \int_{ \mathbb{R} } \len \left| V \right|^{ 1/2 } e^{ itA} g
\rin^2 dt = \int | q ( x ) | \, \left|  g ( x - t ) \right|^2  dx
\,  dt = \len g \rin^2 \int | q ( x ) | dx = \infty \] if $ g \ne
0 $. This proves \textit{ii}). The assertion \textit{iii}) is a
corollary of the following result from \cite{BS},

{\it For any $ \delta $, $ 1 < \delta < 2 $, and any functions $ f
, g $ satisfying \[ \sum_n \( \int_n^{ n + 1 } | f |^2  \)^{
\delta / 2 } < \infty, \; \sum_n \( \int_n^{ n + 1 } | f |^2  \)^{
\delta / 2 } < \infty, \] the operator $ T $ in $ H $ defined by
\[ (Tu) ( x )  = \int_\mathbb{R}  f ( x ) e^{ ixy } g ( y ) u ( y
) dy \] belongs to $ \mathbf{S}^\delta $.}

Applied to $ f = q $ and $ g ( y ) = \( y - z \)^{ -1 } $, $ \im z
\ne 0 $, this result shows that the operator $ V \( A - z \)^{ -1
} \in \mathbf{S}^p $, $ p > 1 $, which implies \textit{iii}).
\end{proof}

\bigskip

\S 2. Let $ T $ be a bounded operator such that $ \sigma ( T )
\subset \mathbb{T} $.

\begin{definition}
The closure of the linear set of vectors $ u \in H $ such that for
all $ v \in H $ the nontangential limits
\[ f_{ u ,v }^\pm ( z ) = \mathop{\lim_{ w \to z }}_{ \left| w \right|^{\pm 1}
\in \mathbb{D} } \llangle \( T - w \)^{ -1 } u , v \rrangle
\] exist and coincide for a.e. $ z \in \mathbb{T} $,
is called the singular subspace of the ope\-rator $ T $. It is
denoted by $ H_s ( T ) $.
\end{definition}

As discussed in the Introduction, the duality problem
\cite{Tikhang} is the question whether the equality \be\la{dua} \(
H_{ac}^w (T) \)^\perp = H_s ( T^* ) \ee holds.

\medskip

{{\parindent=0cm \bf Proposition \cite[Proposition 6.7]{MV}.} {\it
If $ I - T^* T \in \mathbf{S}^1 $, then (\ref{dua}) is satisfied
\footnote{In fact, the cited proposition in \cite{MV} establishes
(\ref{dua}) with the strong a. c. subspace in the place of $
H_{ac}^w (T) $. The latter and the former coincide when $ I - T^*
T \in \mathbf{S}^1 $.}.}

\medskip

More generally, (\ref{dua}) is known to hold if the characteristic
function has weak boundary values a.e. \cite{Ryzov}. We are now
going to construct an example where this property fails. In the
notation of definition \ref{dewcon} let
\[ CN ( T ) \stackrel{ \mbox\small{def}}{ =
} \mathop{Clos} \( \wt{H_+^w} ( T ) \vee \wt{H_-^w} ( T )\) .
\]

Let $ \{ \rho_n \} $, $ n \ge 0 $, be a monotone sequence of
positive numbers tending to $ 0 $, and $ R $ be an operator in $
l^2 ( \Z ) $ defined by $ R e_n = \rho_{ |n| }  e_n $, $ n \in \Z
$. Let $ H = l^2 ( \Z ) \oplus l^2 ( \Z ) $, and let $ U $ be the
operator of right shift in $ l^2 ( \Z ) $, $ U e_n = e_{ n + 1 }
$. Define an operator $ T $ in $ H $ by
\[ T =
\begin{pmatrix} U & R \cr 0 & U \end{pmatrix}. \] Obviously, $ \sigma ( T ) = \mathbb{T}
$.

\begin{theorem} \la{dualexample}
Let $ \{ \rho_j \} \notin l^1 $. Then the operator $ T $ obeys the
following condi\-tions,

\begin{itemize}
\item[\textit{i})] $ CN ( T ) \ne H $,
\item[\textit{ii})] $ H_s ( T^* ) = \{ 0 \} $,
\item[\textit{iii})] $ T = T_0 + S $, where $ T_0 $ is unitary
and $ S $ is an operator whose singular numbers, $ \mu_n ( S ) $,
satisfy $ \mu_n ( S ) \le \rho_{ [ n/2 ] } $.
\end{itemize}
\end{theorem}

\begin{proof} The assertion (\textit{iii}) is obvious
(in fact, $ \mu_{2n} ( S ) = \mu_{ 2n + 1 } ( S ) =\rho_n $ for $
n \ge 1 $). Then, an ele\-men\-ta\-ry calculation gives for any $
f , g \in H $ and $ \lambda \notin \mathbb{T} $, \bequnan \llangle
\( T^* - \lambda \)^{ -1 } f , g \rrangle & = & \llangle \( U^* -
\lambda \)^{ -1 } f_1 , g_1 \rrangle +
\\ & & \quad \llangle \( U^* - \lambda \)^{ -1 } f_2 - \( U^* - \lambda
\)^{ -1 } R \( U^* - \lambda \)^{ -1 } f_1 , g_2 \rrangle .
\eequnan Here
\[ f =
\begin{pmatrix} f_1 \cr f_2
\end{pmatrix} ; \; g
=
\begin{pmatrix} g_1 \cr g_2 \end{pmatrix} . \]

Let $ f $ be from the dense set in the definition of $ H_s ( T^* )
$. Considering the $ g $'s with $ g_2 = 0 $ and arbitrary $ g_1 $
we conclude that $ f_1 = 0 $, since $ U $ is an absolutely
continuous unitary operator. Then, by the same reason, considering
the $ g $'s with $ g_1 = 0 $ we obtain that $ f_2 = 0 $, that is,
$ H_s ( T^* ) $ is trivial. Also, the absolute continuity of $ U $
implies that $ CN ( T ) \supset  { l^2 ( \mathbb{Z} ) \choose 0 }
$. Let us show that in fact $ CN ( T ) = { l^2 ( \mathbb{Z} )
\choose 0 } $. Actually, we shall show that if for a $ u = { u_1
\choose u_2 } \in H $ the function \be\la{normu} \llangle \( T -
\lambda \)^{ -1 } u ,
\begin{pmatrix} e_j \cr 0 \end{pmatrix} \rrangle =
\llangle \( U - \lambda \)^{ -1 } u_1 , e_j \rrangle  - \llangle
\( U - \lambda \)^{ -1 } R \( U - \lambda \)^{ -1 } u_2 , e_j
\rrangle \ee is in $ H^2 $ for all $ j $, and the $ H^2 $-norm of
it is bounded above in $ j $, then $ u_2 = 0 $. Since for a $ u\in
\wt{H_+^w} $ this norm is indeed bounded in $ j $, this is going
to imply that $ \wt{H_+^w} \subset { l^2 ( \mathbb{Z} ) \choose 0
} $.

Taking into account that
\[ \len \llangle \( U - \cdot \)^{ -1 } u_1 , e_j \rrangle \rin_{
H^2 }^2 = \sum_1^\infty \left| u_{1 , k + j} \right|^2 \le \len u
\rin^2 ,
\] we only have to check that the $ H^2 $-norm of the
second term in the right hand side in (\ref{normu}) is unbounded
as $ j \to \infty $ if $ u_2 \ne 0 $. Indeed, a straightforward
calculation gives that for $ \lambda \in \mathbb{D} $
\[ \llangle \( U - \lambda \)^{ -1 } R \( U - \lambda \)^{ -1 }
u_2 , e_j \rrangle = \sum_{ s = 0 }^\infty \lambda^s u_{2, s+j+2 }
\sum_{ m =1 }^{ s+1 } \rho_{ | m+j | } . \] Thus, the $ H^2 $
-norm of the left hand side is \[ \sum_{ s = 0 }^\infty \left|
u_{2, s+j+2 } \right|^2 \left| \sum_{ m = j + 1 }^{ j+s+1 } \rho_{
|m| } \right|^2 .
\] Suppose that $ u_{ 2 , r } \ne 0 $ for some $ r $. Then for $ j
$ negative enough this norm is bounded below by $ \left| u_{2, r }
\right|^2 \left| \sum_{ m =j+1 }^{ r-1 } \rho_{ |m| } \right|^2 $
which goes to infinity when $ j \to - \infty $ by the assumption
about $ \rho_j $. The inclusion $ \wt{H_-^w} \subset { l^2 (
\mathbb{Z} ) \choose 0 } $ is checked similarly.
\end{proof}

In the example constructed the subspaces $ CN ( T ) $ and $ H_{
ac}^w ( T ) $ coincide, thus showing that even the condition $ CN
( T )^\perp \subset H_s ( T^* ) $, weaker than (\ref{dua}) in
general, can fail for non-trace class perturbations. It would be
interesting to see if the linear resolvent growth condition
(\ref{lrg}), violated in the example above, implies the inclusion
$ CN ( T )^\perp \subset H_s ( T^* ) $.

\bigskip

\S 3. We now turn to the case $ p \ne 2 $.

\begin{lemma} \la{incl}
$ \overline{\( L - z_0 \)^{ -1 } H_{ac}^{w, p_1  }} = H_{ac}^{w,
p_1  } \subset H_{ac}^{w, { p_2 } } $ for all $ z_0 \in \rho ( L )
$ and $ 1 \le p_2 \le p_1 \le 2 $ save for $ p_1 = 1 $.
\end{lemma}

\begin{proof}
An application of the H\"{o}lder inequality to the resolvent
identity \[ f_{ \( L - z_0 \)^{ -1 } u , v } ( z ) = \frac 1{ z -
z_0 } \( f_{ u , v } ( z )  -   f_{ u , v } ( z_0 ) \)  \] shows
that \[ \( L - z_0 \)^{ -1 } \wt{ H_{ac}^{w, p_1 } } \subset \wt{
H_{ac}^{w, { p_2 } }}  \] for all $ z_0 \in \rho ( L ) $. Since
the linear set in the left hand side is independent of the choice
of $ z_0 \in \rho ( L ) $, the asymptotics (\ref{inft}) implies
that the set is dense in $ H_{ac}^{w, p_1 } $.
\end{proof}

It is also easy to check that \be\la{estsm} \len \( L - z \)^{ -1
} u \rin \le C_u \left| \im z \right|^{ -1/p } \ee for any $ u \in
\wt{ H_{ac}^{w, p } } $, $ 1 < p \le 2 $.

\begin{theorem}
\la{acacwHp}
If $ L $ is dissipative ($ V \ge 0 $) and $ 1 < p \le 2 $, then $
H_{ac}^{w, p } ( L ) = H_{ac}^w ( L ) $.
\end{theorem}

In the earlier note \cite{equival} we proved that if $ L $ is
dissipative then $ H_{ac} ( L ) = H_{ac}^w ( L ) $.

\begin{proof}
Since the inclusion $ H_{ac}^w \subset H_{ac}^{w, p } $ is
contained in lemma \ref{incl}, it remains to check that $
H_{ac}^{w, p } \subset H_{ ac } $. We do so by proving that any
vector $ u $ from a dense subset in $ H_{ac}^{w, p } $ is a strong
smooth vector. Recall \cite{N} that for a dissipative operator $ L
$ the restriction of the function $ V^{ 1/2 } \( L - z \)^{ -1 } w
$ to $ \mathbb{C}_- $ belongs to $ \H^2_- $ for all $ w \in H $.
Hence, we only have to verify that $ \( V^{ 1/2 } \( L - z \)^{ -1
} u \)_+ \in \H^2_+ $ for all $ u $ from the dense subset.

Let $ u \in \( L + i \)^{ -2 } \wt{H_{ac }^{w,p}} $. Taking into
account (\ref{estsm}) one easily checks that the function $ \( L -
\cdot - i \von \)^{ -1 } u \in L^2 ( \mathbb{R}, H ) $ for any $
\von > 0 $. We are first going to show that \be\la{epsq} \sup_{
\von > 0 } \( \von \int_\mathbb{R} \len \( L - k - i \von \)^{ -1
} u \rin^2 dk \) < \infty . \ee

For $ \von > 0 $ and $ t < 0 $ we define \be\la{ut} u ( t ) = -
\frac 1{ 2 \pi i } \mathop{\lim}_{ N \to \infty } \int_{ -N }^N
e^{ i ( k + i \von ) t } \( L - k - i \von \)^{ -1 } u \, dk . \ee
Then

$ 1^\circ $. The limit in (\ref{ut}) exists for all $ t < 0 $, and
is independent of $ \von > 0 $.

$ 2^\circ $. $ \sup_{ t < 0 } \| u ( t ) \| < \infty $.

The assertion $ 1^\circ $ follows immediately from the possibility
to rewrite (\ref{ut}) in the form ($ \lambda = k + i \von $)

\[ u (t) = - \frac 1{ 2 \pi i }
\int_{ \mathbb{R} } \frac{e^{ i \lambda t }}{ \( \lambda + i \)^2
} \( L - \lambda \)^{ -1 } u_2 \, dk - e^t ( it \, u_1 + u ) \]
where $ u_1 = ( L + i ) u $, $ u_2 = ( L + i )^2 u $.

Let us establish $ 2^\circ $. For any $ v \in H $ the scalar
product $ \llangle u ( t ) , v \rrangle $ equals to \[ - \frac 1{
2 \pi i } \int_{ \mathbb{R} } \frac{e^{ i \lambda t }}{ \( \lambda
+ i \)^2 } f_{ u_2 , v } ( \lambda ) \, dk + r_t \] where $ | r_t
| \le C  \| v \| $, the bound being uniform in $ t < 0 $. Since $
f_{ u_2 , v } \in H^p $, one can pass to the limit $ \von \to 0 $
in this integral, and use the H\"{o}lder inequality to estimate
the modulus of it by
\[ C \len f_{ u_2 , v } \rin_{ H^p } \le C \| v \| \] with a
constant $ C $ independent of $ t $ and $ v $. This implies $
2^\circ $.

Let $ \cal F $ stands for the conjugate Fourier transform in $ L^2
( \mathbb{R}, H ) $. Then (\ref{ut}) means that the restriction of
the function
\[ \Psi ( t ) = - \frac 1{
i \sqrt{ 2 \pi} } \, {\cal F} \left[ \( L - \cdot - i \von \)^{ -1
} u \right]
\] to $ t < 0 $ coincides with $ e^{ \von t } u ( t ) $. On the other
hand, $ \Psi ( t ) = 0 $ for $ t > 0 $ by the Paley - Wiener
theorem. Applying the Parseval equality and taking into account
the property $ 2^\circ $, we find
\[ \int_\mathbb{R} \len \( L - k - i \von \)^{ -1 } u \rin^2 dk = C
\int_{ - \infty }^0 e^{ 2\von t } \len u ( t ) \rin^2 \, dt \le C
\von^{ -1 } . \] The estimate (\ref{epsq}) is proved.

We are now going to use the following easily verified identity
valid for all $ \lambda \in \rho ( L ) $ ($ \von = \im \lambda $),
\[ \len V^{ 1/2 } \( L - \lambda \)^{ -1 } u \rin^2 = \von \len \( L
- \lambda \)^{ -1 } u \rin^2 - \im f_{ u , u } ( \lambda ) . \]
The second term in the right hand side can be rewritten in the
form \[ \im \left[ \frac 1{ \lambda +i } \( f_{ u_1 , u } (
\lambda ) - f_{ u_1 , u } ( -i ) \) \right] \] which implies that
the term is conditionally integrable in $ \re \lambda $ over the
real line, and the integrals are uniformly bounded in $ \von $.
Together with (\ref{epsq}), this shows that $ u \in \wt{H_{ac}} $.
It remains to notice that $\( L + i \)^{ -2 } \wt{H_{ac }^{w,p}} $
is dense in $ H_{ac }^{w,p} $.
\end{proof}

\begin{remark} \la{ple1} In a similar way, the subspaces $ H_{ac}^{w, p } ( T ) $
can be defined for pertur\-ba\-tions of unitary operators. These
subspaces can be shown to coincide with $ H_{ac} ( T ) $ for all $
p $, $ 1 \le p \le 2 $. The subspaces $ H_{ac}^{w, p } ( T ) $ can
also be defined for $ 0 < p < 1 $ but they coincide with $ H $ for
any unitary operator $ T $, rendering the definition meaningless.
The reason is that the Cauchy transform of any finite measure is
in $ H^p ( \mathbb{D} ) $ for $ 0 < p < 1 $. In the real line
context, the definition of the subspace for $ p = 1 $ requires a
regularization at infinity.
\end{remark}

\section*{Appendix}

\begin{proposition}
Let $ L $ be a selfadjoint operator and let $\mathcal{H}_{ac} ( L
) $ be its a. c. subspace defined via the spectral theorem. Then $
H_{ac}^{w, p } ( L ) = \mathcal{H}_{ac} ( L ) $ for all $ p \in (
1 , 2 ] $.
\end{proposition}

\begin{proof}
Let $ d\mu_{u,v} ( t ) $ be the matrix element of the spectral
measure of $ L $ on vectors $ u , v \in H $. Then $$ f_{ u , v } (
z ) = \int_\R \frac 1{ t - z } \, d\mu_{u,v} ( t ) $$ for all $ u
, v \in H $. Let $ u = \( L - z_0 \)^{ -1 } w $ with a $ w \in
\wt{H_{ac}^{w, p }} $ and $ z_0 \in \rho ( L ) $, then $ f_{ u , v
} $ is represented as the Cauchy transform of its boundary values,
$$ f_{ u , v } ( z ) = \int_\R \frac 1{ (t - z) (t-z_0) } \,
f_{w,v} ( t ) dt . $$ Notice that $ \( t-z_0 \)^{ -1 } f_{w,v} ( t
) \, dt $ is a finite Borel mea\-su\-re. Com\-paring the two
representations and using the Riesz bro\-thers theorem, one
concludes that the measure $$ d \mu_{u,v} - \( t-z_0 \)^{ -1 }
f_{w,v} ( t ) dt $$ is a. c. for all $ v $. Hence $ d \mu_{u,v} $
is a. c. as well. Since the set of such $ u $'s is dense in $
H_{ac}^{w, p } $, the inclusion $ H_{ac}^{w, p } \subset
\mathcal{H}_{ac} $ follows. The inclusion $ u \in H_{ac}^{w, p } $
is obvious for any $ u \in \mathcal{H}_{ac} $ satisfying $
d\mu_{u,u}/ dt \in L^\infty ( \mathbb{R} ) $. Since the set of
such $ u $'s is dense in $ u \in \mathcal{H}_{ac} $, it implies
that $ \mathcal{H}_{ac} \subset H_{ac}^{w, p } $.
\end{proof}

\bigskip
{\bf Acknowledgements.} The author is indebted to S. Naboko for
helpful discussions, and to an unknown referee for remarks and
suggestions. The assertion of theorem \ref{acacwHp} is the answer
to a question of E. Shargorodsky. The work is partially supported
by the INTAS Grant 05-1000008-7883 and the RFBR Grant 06-01-00249.

\end{document}